\documentclass[12pt]{article}
\usepackage{graphicx}

\newcommand{\hide}[1]{}

\newcommand{\ac}[1]{{\tt #1}}  
\newcommand{\inst}[1]{{\tt #1}}  

\newcommand{\davg}{\Delta_{\mbox{\tiny\em avg}}}

\newcommand{\dlbs}{don't look bits}

\newcommand{\twoopt}{{\tt 2-opt}}
\newcommand{\threeopt}{{\tt 3-opt}}
\newcommand{\lk}{{\tt LK}}
\newcommand{\kl}{{\tt KL}}

\newcommand{\pscProc}{{\bf procedure}}
\newcommand{\pscEnd}{{\bf end}}
\newcommand{\pscRepeat}{{\bf repeat}}
\newcommand{\pscUntil}{{\bf until}}

\newcommand{\genls}{{\small\sf LocalSearch}}
\newcommand{\genps}{{\small\sf Per\-tur\-ba\-tion}}
\newcommand{\genkm}{{\small\sf Per\-tur\-ba\-tion}}
\newcommand{\genac}{{\small\sf AcceptanceCriterion}}
\newcommand{\genis}{{\small\sf GenerateInitialSolution}}

\newcommand{\ls}{{local search}}

\newcommand{\fsp}{{flow shop problem}}
\newcommand{\afsp}{{FSP}}
\newcommand{\jsp}{{job shop scheduling problem}}
\newcommand{\ajsp}{{JSP}}
\newcommand{\qap}{{quadratic assignment problem}}
\newcommand{\tsp}{{traveling salesman problem}}

\newcommand{\ils}{{iterated local search}}
\newcommand{\ails}{{ILS}}
\newcommand{\vns}{{variable neighborhood search}}
\newcommand{\avns}{{VNS}}

\newcommand{\sa}{{simulated annealing}}
\newcommand{\asa}{{SA}}
\newcommand{\ts}{{tabu search}}
\newcommand{\ats}{{TS}}

\newcommand{\ma}{{memetic algorithms}}

\newcommand{\agrasp}{{GRASP}}
\newcommand{\aco}{{ant colony optimization}}
\newcommand{\aaco}{{ACO}}

\begin{document}

\title{Iterated Local Search}

\author{Helena R. Louren\c co\\
Universitat Pompeu Fabra, Barcelona, Spain\\
helena.ramalhinho@econ.upf.es\\
Olivier C. Martin\\
Universit\'e Paris-Sud, Orsay, France\\
martino@ipno.in2p3.fr\\
Thomas St\"utzle\\
Darmstadt University of Technology, Darmstadt, Germany\\
stuetzle@informatik.tu-darmstadt.de\\
}

\bigskip
\date{{\it To appear in ``Handbook on MetaHeuristics'', \\
Ed. F. Glover and G. Kochenberger}}

\maketitle

\section{Introduction}
\label{s:i}

The importance of high performance algorithms for tackling difficult
optimization problems cannot be understated, and in many cases the
only available methods are metaheuristics. When designing a
metaheuristic, it is preferable that it be simple, both conceptually
and in practice. Naturally, it also must be effective, and if
possible, general purpose.
If we think of a metaheuristic as simply a construction for guiding
(problem-specific) heuristics, the ideal case is when the
metaheuristic can be used without {\it any} problem-dependent
knowledge.

As metaheuristics have become more and more sophisticated, this ideal
case has been pushed aside in the quest for greater performance. As a
consequence, problem-specific knowledge (in addition to that built
into the heuristic being guided) must now be incorporated into
metaheuristics in order to reach the state of the art
level. Unfortunately, this makes the boundary between heuristics and
{\em meta\/}heuristics fuzzy, and we run the risk of loosing both
simplicity and generality. To counter this, we appeal to modularity
and try to decompose a metaheuristic algorithm into a few parts, each
with its own specificity. In particular, we would like to have a
totally general purpose part, while any problem-specific knowledge
built into the metaheuristic would be restricted to another
part. Finally, to the extent possible, we prefer to leave untouched
the embedded heuristic (which is to be ``guided'') because of its
potential complexity. One can also consider the case where this
heuristic is only available through an object module, the source code
being proprietary; it is then necessary to be able to treat it as a
``black-box'' routine. Iterated local search provides a simple way to
satisfy all these requirements.

The essence of the \ils\ metaheuristic can be given in a nut-shell:
one {\it iteratively} builds a sequence of solutions generated by the
embedded heuristic, leading to far better solutions than if one were
to use repeated random trials of that heuristic.  This simple
idea~\cite{Baxter81} has a long history, and its rediscovery by many
authors has lead to many different names for \ils\ like {\it iterated
descent}~\cite{Baum_86a,Baum_86b}, {\it large-step Markov
chains}~\cite{MarOttFel91:cs}, {\it iterated
Lin-Kernighan}~\cite{Joh90}, {\it chained local
optimization}~\cite{MarOtt96:aor}, or combinations of
these~\cite{AppCooRoh99}~...  Readers interested in these historical
developments should consult the review~\cite{JohMcG97}.  For us, there
are two main points that make an algorithm an \ils: (i) there must be
a single chain that is being followed (this then excludes
population-based algorithms); (ii) the search for better solutions
occurs in a reduced space defined by the output of a black-box
heuristic. In practice, local search has been the most frequently used
embedded heuristic, but in fact any optimizer can be used, be-it
deterministic or not.

The purpose of this review is to give a detailed description of
\ils\ and to show where it stands in terms of performance. So far, in
spite of its conceptual simplicity, it has lead to a number of
state-of-the-art results without the use of too much problem-specific
knowledge; perhaps this is because \ils\ is very malleable, many
implementation choices being left to the developer.  We have
organized this chapter as follows. First we give a high-level
presentation of \ils\ in Section~\ref{sect_iterating}  Then we
discuss the importance of the different parts of the metaheuristic in
Section~\ref{sect_making}, especially the subtleties associated with
perturbing the solutions.  In Section~\ref{s:rw} we go over past work
testing \ils\ in practice, while in Section~\ref{s:r} we discuss
similarities and differences between \ils\ and other metaheuristics.
The chapter closes with a summary of what has been achieved so far and
an outlook on what the near future may look like.

\section{Iterating a local search}
\label{sect_iterating}

\subsection{General framework}
\label{subsect_framework}

We assume we have been given a problem-specific 
heuristic optimization algorithm 
that from
now on we shall refer to as a \ls\ (even if in fact it is not
a true \ls). This algorithm is implemented via a computer 
routine that we call
{\genls}. The question we ask is ``Can such an algorithm be improved
by the use of iteration?''. Our answer is ``YES'', and in fact the
improvements obtained in practice are usually significant. Only in
rather pathological cases where the iteration method is
``incompatible'' with the local search will the improvement be minimal. In
the same vein, in order to have the {\it most} improvement possible,
it is necessary to have some understanding of the way the
\genls\ works.  However, to keep this presentation as simple as
possible, we shall ignore for the time being these complications; the
additional subtleties associated with tuning the iteration to the
\ls\ procedure will be discussed in Section~\ref{sect_making}
Furthermore, all issues associated with the actual speed of the
algorithm are omitted in this first section as we wish to focus solely
on the high-level architecture of \ils .

Let $\cal C$ be the cost function of our combinatorial optimization
problem; $\cal C$ is to be {\it minimized}. We label candidate
solutions or simply ``solutions'' by $s$, and denote by $\mathcal S$
the set of all $s$ (for simplicity $S$ is taken to be finite, but it
does not matter much).  Finally, for the purposes of this high-level
presentation, it is simplest to assume that the local search procedure
is deterministic and memoriless:\footnote{The reader can check that
very little of what we say really uses this property, and in practice,
many successful implementations of \ils\ have non-deterministic local
searches or include memory.}  for a given input $s$, it always outputs
the same solution $s^*$ whose cost is less or equal to ${\cal
C}(s)$. \genls\ then defines a many to one mapping from the set $\cal
S$ to the smaller set ${\mathcal S}^*$ of locally optimal solutions
$s^*$. To have a pictorial view of this, introduce the ``basin of
attraction'' of a local minimum $s^*$ as the set of $s$ that are
mapped to $s^*$ under the local search routine. \genls\ then takes one
from a starting solution to a solution at the bottom of the
corresponding basin of attraction.

\begin{figure}[tb]
\begin{center}
	\includegraphics[width=0.75\textwidth]{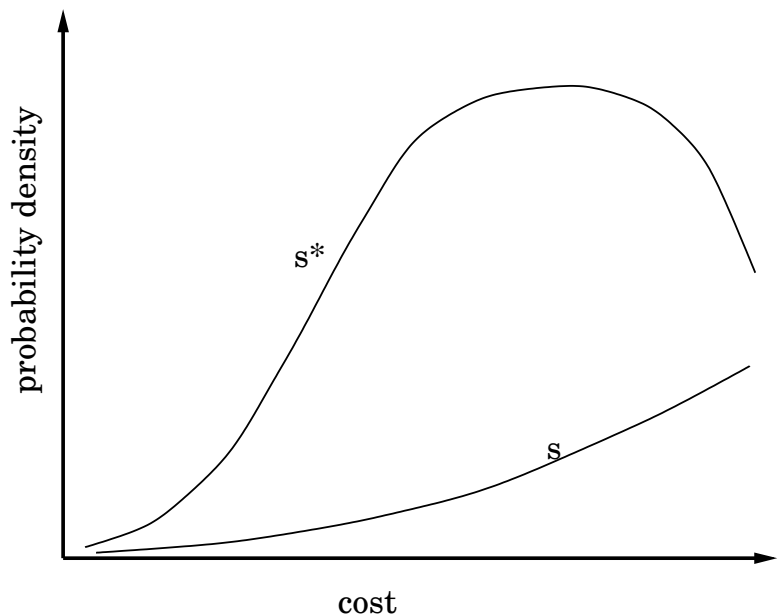}
\end{center}
\caption{Probability densities of costs. The curve labeled
$s$ gives the cost density for all solutions, while the curve
labeled $s^*$ gives it for the solutions that are local optima.}
\label{fig_costs}
\end{figure}

Now take an $s$ or an $s^*$ at random. Typically, the distribution of
costs found has a very rapidly rising part at the lowest values.  In
Figure~\ref{fig_costs} we show the kind of distributions found in
practice for combinatorial optimization problems having a finite
solution space. The distribution of costs is bell-shaped, with a mean
and variance that is significantly smaller for solutions in ${\cal
S}^*$ than for those in $\mathcal S$. As a consequence, it is much better
to use \ls\ than to sample randomly in $\mathcal S$ if one seeks low cost
solutions. The essential ingredient necessary for \ls\ is a
neighborhood structure. This means that $\mathcal S$ is a ``space'' with
some topological structure, not just a set.  Having such a space
allows one to move from one solution $s$ to a better one in an
intelligent way, something that would not be possible if $\mathcal S$ were just
a set.

Now the question is how to go beyond this use of \genls . More precisely,
given the mapping from $\mathcal S$ to ${\cal S}^*$, how can one further
reduce the costs found without opening up and modifying
\genls , leaving it as a ``black box'' routine?

\subsection{Random restart}
\label{subsect_multi-start}
The simplest possibility to improve upon a cost found by \genls\ is to
repeat the search from another starting point.  Every $s^*$ generated
is then independent, and the use of multiple trials allows one to
reach into the lower part of the distribution. 
Although such a ``random restart''
approach with independent samplings is sometimes a useful strategy (in
particular when all other options fail), it breaks down as the
instance size grows because in that limit the tail of the distribution
of costs collapses.  Indeed, empirical studies~
\cite{JohMcG97} and general arguments~\cite{SM99} indicate that
local search algorithms on large generic instances lead to
costs that: (i) have a mean that is a fixed
percentage excess above the optimum cost; (ii) have a {\it
distribution} that becomes arbitrarily peaked about the
mean when the instance size goes to 
infinity. This second property makes it impossible in practice to find an
$s^*$ whose cost is even a little bit lower percentage-wise 
than the typical cost.
Note however that there do exist many solutions of significantly lower cost,
it is just that {\it random} sampling has a lower and lower
probability of finding them as the instance size increases. To reach
those configurations, a biased sampling is necessary; this is
precisely what is accomplished by a stochastic search.

\subsection{Searching in ${\cal S}^*$}
\label{subsect_searching}
To overcome the problem just mentioned associated with large instance
sizes, reconsider what local search does: it takes one from $\mathcal S$
where $\cal C$ has a large mean to ${\mathcal S}^*$ where $\cal C$ has a
smaller mean. It is then most natural to invoke recursion: use \ls\ to
go from ${\mathcal S}^*$ to a smaller space ${\cal S}^{**}$ where the mean
cost will be still lower! That would correspond to an algorithm with
one local search nested inside another. Such a construction could
be iterated to as many levels as desired, leading to a hierarchy of
nested local searches. But upon closer scrutiny, we see that the
problem is precisely how to formulate local search beyond the lowest
level of the hierarchy: local search requires a neighborhood structure
and this is not {\it \`a priori} given.  The fundamental difficulty is to
define neighbors in ${\mathcal S}^*$ so that they can be enumerated and
accessed efficiently. Furthermore, it is desirable for nearest
neighbors in ${\mathcal S}^*$ to be relatively close when using the
distance in $\mathcal S$; if this were not the case, a stochastic search
on ${\mathcal S}^*$ would have little chance of being effective.

Upon further thought, it transpires that one can introduce
a good neighborhood structure on ${\mathcal S}^*$ as follows. First,
one recalls that a neighborhood structure on a set
$\mathcal S$ directly induces a neighborhood structure on {\it subsets} of
$\mathcal S$: two subsets are nearest neighbors simply if they contain
solutions that are nearest neighbors. Second, take these subsets to be
the basins of attraction of the $s^*$; in effect, we are lead to identify any
$s^*$ with its basin of attraction. This then immediately gives the
``canonical'' notion of neighborhood on ${\mathcal S}^*$, notion which can be
stated in a simple way as follows: $s_1^*$ and $s_2^*$ are neighbors in
${\mathcal S}^*$ if their basins of attraction ``touch'' ({\it i.e.},
contain nearest-neighbor solutions in $\mathcal S$). Unfortunately this
definition has the major drawback that one cannot in practice list the
neighbors of an $s^*$ because there is no computationally efficient
method for finding all solutions $s$ in the basin of attraction of
$s^*$. Nevertheless, we can {\it stochastically} generate nearest
neighbors as follows.  Starting from $s^*$, create a randomized path
in $\mathcal S$, $s_1$, $s_2$, ..., $s_i$, where $s_{j+1}$ is a nearest
neighbor of $s_j$. Determine the first $s_j$ in this path that belongs
to a different basin of attraction so that applying \ls\ to $s_j$
leads to an $s^*{'} \ne s^*$. Then $s^*{'}$ is a nearest-neighbor of
$s^*$.

Given this procedure, we can in principle
perform a local search\footnote{Note that the local
search finds neighbors stochastically; generally there is no efficient
way to ensure that one has tested {\it all} the neighbors of any given
$s^*$.} in ${\mathcal S}^*$. Extending the
argument recursively, we see that it would be possible to have an
algorithm implementing nested searches, performing local search on
$\mathcal S$, ${\cal S}^*$, ${\cal S}^{**}$, etc... in a hierarchical way.
Unfortunately, the implementation of nearest neighbor search at the
level of ${\mathcal S}^*$ is much too costly computationally because of
the number of times one has to execute \genls.
Thus we are led to abandon the (stochastic)
search for {\it nearest neighbors} in ${\mathcal S}^*$; instead we use a
weaker notion of closeness which then allows for a fast stochastic
search in ${\mathcal S}^*$.
Our construction leads to a (biased) sampling of ${\mathcal S}^*$; such a
sampling will be better than a random one if it is possible to find
appropriate computational ways to go from one $s^*$ to another.
Finally, one last advantage of this modified notion of closeness
is that it does not require basins of attraction to be defined;
the local search can then incorporate memory or be non-deterministic,
making the method far more general.

\subsection{Iterated Local Search}
\label{subsect_ILS}
We want to explore ${\mathcal S}^*$ using a walk that steps from one $s^*$
to a ``nearby'' one, without the constraint of using only nearest
neighbors as defined above. Iterated local search (ILS) achieves this
heuristically as follows. Given the current $s^*$, we first apply a
change or perturbation that leads to an intermediate state $s'$ (which
belongs to $\mathcal S$). Then \genls\ is applied to 
$s{'}$ and we reach a solution
$s^*{'}$ in ${\mathcal S}^*$. If $s^*{'}$ passes an acceptance test, it
becomes the next element of the walk in ${\mathcal S}^*$; otherwise, one
returns to $s^*$.  The resulting walk is a case
of a stochastic search in ${\mathcal S}^*$, but where neighborhoods
are never explicitly introduced. This \ils\ procedure should lead to
good biased sampling as long as the perturbations are neither too
small nor too large.  If they are too small, one will often
fall back to $s^*$ and few new solutions of ${\mathcal S}^*$
will be explored.  If on the contrary the perturbations are too large,
$s{'}$ will be random, there will be no
bias in the sampling, and we will recover a random restart type
algorithm.

\begin{figure}[tb]
\centerline{
	\includegraphics[width=0.75\textwidth]{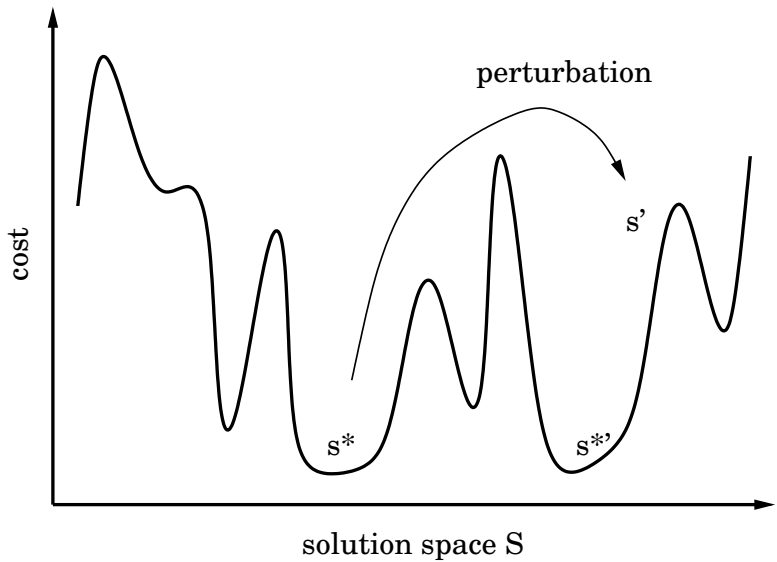}}
\caption{Pictorial representation of \ils . Starting with a local minimum
$s^*$, we apply a perturbation leading to a solution
$s{'}$. After applying \genls , we find a new local
minimum $s^*{'}$ that may be better than $s^*$.}
\label{fig_landscape}
\end{figure}

The overall ILS procedure is pictorially illustrated in
Figure~\ref{fig_landscape}.  To be complete, let us note that
generally the \ils\ walk will not be reversible; in particular one may
sometimes be able to step from $s_1^*$ to $s_2^*$ but not from $s_2^*$ to
$s_1^*$.  However this ``unfortunate'' aspect of the procedure does
not prevent ILS from being very effective in practice.

Since deterministic perturbations may lead to short cycles (for
instance of length $2$), one should randomize the perturbations or
have them be adaptive so as to avoid this kind of cycling. If the
perturbations depend on any of the previous $s^*$, one has a walk in
${\mathcal S}^*$ with {\it memory}. Now the reader may have noticed that
aside from the issue of perturbations (which use the structure on
$\mathcal S$), our formalism reduces the problem to that of a stochastic
search on ${\mathcal S}^*$. Then all of the bells and whistles
(diversification, intensification, tabu, adaptive perturbations and
acceptance criteria, etc...)  that are commonly used in that context
may be applied here.  This leads us to define \ils\ algorithms as
metaheuristics having the following high level architecture:

\begin{quote}
\begin{tabbing}
----\=----\=----\=----\=----\=----\=----\=----\=\kill
\pscProc\ {\it Iterated Local Search} \\
\> $s_0 = \mbox{\sf GenerateInitialSolution}$\\
\> $s^* = \mbox{\sf LocalSearch}(s_0)$\\
\> \pscRepeat         \\
\> \> $s' = \mbox{\sf Perturbation}(s^*,\mbox{\em history\/})$\\
\> \> $s^*{'} = \mbox{\sf LocalSearch}(s')$\\
\> \> $s^* = \mbox{\sf\ AcceptanceCriterion}
(s^*,s^*{'},\mbox{\em history\/})$\\
\> \pscUntil\ termination condition met\\
\pscEnd \\
\end{tabbing}
\label{f:ils} 
\end{quote}

In practice, much of the potential complexity of ILS is hidden in
the history dependence. If there happens to be no such dependence, the
walk has no memory:\footnote{Recall that to simplify this
section's presentation, the local search is assumed to have
no memory.}
 the perturbation and acceptance criterion do not
depend on any of the solutions visited previously during the walk, and
one accepts or not $s^*{'}$ with a fixed rule. This leads to random
walk dynamics on ${\mathcal S}^*$ that are ``Markovian'', the probability
of making a particular step from $s_1^*$ to $s_2^*$ depending only on
$s_1^*$ and $s_2^*$. Most of the work using ILS has been of this
type, though recent studies show unambiguously that incorporating
memory enhances performance~\cite{Stu98:phd}.

Staying within Markovian walks, the most basic acceptance criteria
will use only the difference in the costs of $s^*$ and $s^*{'}$; this
type of dynamics for the walk is then very similar in spirit to what
occurs in simulated
annealing. A limiting case of this is to accept only improving moves,
as happens in simulated annealing at zero temperature; the algorithm
then does (stochastic) descent in ${\mathcal S}^*$.  If we add to such a
method a CPU time criterion to stop the search for improvements, the
resulting algorithm pretty much has two nested local searches;
to be precise, it has a local search operating on $\mathcal S$ embedded in a
stochastic search operating on ${\mathcal S}^*$. More generally, one can
extend this type of algorithm to more levels of nesting, having a
different stochastic search algorithm for ${\mathcal S}^*$, ${\cal
S}^{**}$ etc...  Each level would be characterized by its own type of
perturbation and stopping rule; to our knowledge, such a construction
has never been attempted.

We can summarize this section by saying that the potential power of
\ils\ lies in its {\it biased} sampling of the set of local
optima. The efficiency of this sampling depends both on the kinds of
perturbations and on the acceptance criteria. Interestingly, even with
the most na\"\i ve implementations of these parts, \ils\ is
much better than random restart.
But still much better
results can be obtained if the \ils\ modules are optimized. First, the
acceptance criteria can be adjusted empirically as in simulated
annealing without knowing anything about the problem being optimized.
This kind of optimization will be familiar to any user of
metaheuristics, though the questions of memory may become quite
complex. Second, the \genps\ routine can incorporate as much
problem-specific information as the developer is willing to put into
it. In practice, a rule of thumb can be used as a guide:
``a good perturbation transforms one excellent solution into an 
excellent starting point for a local search''.
Together, these different aspects show that \ils\ algorithms can have a wide
range of complexity, but complexity may be added progressively and in
a modular way.  (Recall in particular that all of the fine-tuning that
resides in the embedded \ls\ can be ignored if one wants, and it does not
appear in the metaheuristic per-se.) This makes \ils\ an appealing
metaheuristic for both academic and industrial applications.  The
cherry on the cake is speed: as we shall soon see, one can perform $k$
local searches embedded within an \ils\ {\it much} faster than if the
$k$ local searches are run within random restart.

\section{Getting high performance}
\label{sect_making}

Given all these advantages, we hope the reader is now motivated to go
on and consider the more nitty-gritty details that arise
when developing an ILS for a new application. In this section, we
will illustrate the main issues that need to be tackled
when optimizing an ILS in order to achieve high performance.

There are four components to consider: \genis, \genls,
\genkm, and \genac. Before 
attempting to develop a state-of-the-art algorithm, it is relatively
straight-forward to develop a more basic version of ILS. Indeed, (i)
one can start with a random solution or one returned by some greedy
construction heuristic; (ii) for most problems a local search
algorithm is readily available; (iii) for the perturbation, a random
move in a neighborhood of higher order than the one used by the local
search algorithm can be surprisingly effective; and (iv) a reasonable
first guess for the acceptance criterion is to force the cost to
decrease, corresponding to a first-improvement descent in the set
$\mathcal{S}^*$.  Basic ILS implementations of this type usually lead
to much better performance than random restart approaches. The developer
can then run this basic ILS to build his intuition and try to improve
the overall algorithm performance by improving each of the four
modules. This should be particularly effective if it is possible to
take into account the specificities of the combinatorial optimization
problem under consideration. In practice, this tuning is easier for
ILS than for \ma\ or \ts\ to name but these metaheuristics. The reason
may be that the complexity of ILS is reduced by its modularity, the
function of each component being relatively easy to understand.
Finally, the last task to consider is the overall optimization 
of the ILS algorithm; indeed, the different components affect 
one another and so it is necessary to understand their interactions. 
However, because these interactions are so
problem dependent, we wait till the
end of this section before discussing that kind of ``global'' optimization.

Perhaps the main message here is that the developer can choose
the level of optimization he wants. In the absence of any
optimizations, ILS is a simple, easy to implement, and quite effective
metaheuristic. But with further work on its four components, ILS
can often be turned into a very competitive or even
state of the art algorithm.

\subsection{Initial solution}
\label{Sec3:ss:is}

Local search applied to the initial solution $s_0$ gives the starting
point $s_0^*$ of the walk in the set $\mathcal{S}^*$. Starting with
a good $s_0^*$ can be important if high-quality solutions are
to be reached {\it as fast as possible}. 

Standard choices for $s_0$ are either a random initial
solution or a solution returned by a greedy construction
heuristic.
A greedy initial solution $s_0$ has two main advantages over random
starting solutions: (i) when combined with local search, greedy
initial solutions often result in better quality solutions $s_0^*$;
(ii) a local search from greedy solutions takes, on average, less
improvement steps and therefore the local search requires less CPU
time.\footnote{Note that the best possible greedy initial solution
need not be the best choice when combined with a local
search. For example, in \cite{JohMcG97}, it is shown that the
combination of the Clarke-Wright starting tour (one of the
best performing TSP construction heuristics) with local search resulted in
worse local optima than starting from random initial solutions when
using $3$-opt. Additionally, greedy algorithms which generate very
high quality initial solutions can be quite time-consuming.}

The question of an appropriate initial solution for
(random restart) local
search carries over to ILS because of the dependence of the walk in
$\mathcal{S}^*$ on the initial solution $s_0^*$. Indeed, when starting with
a random $s_0$, ILS may take several iterations to catch up 
in quality with runs using an
$s_0^*$ obtained by a greedy initial solution.
Hence, for short computation times the initial solution
is certainly important to achieve the highest quality solutions
possible. For larger computation times, the dependence on
$s_0$ of the final
solution returned by ILS reflects just 
how fast, if at all, the memory of the initial
solution is lost when performing the walk in $\mathcal{S}^*$.

\begin{figure}[tb]
  \begin{center}
	\includegraphics[width=0.5\textwidth]{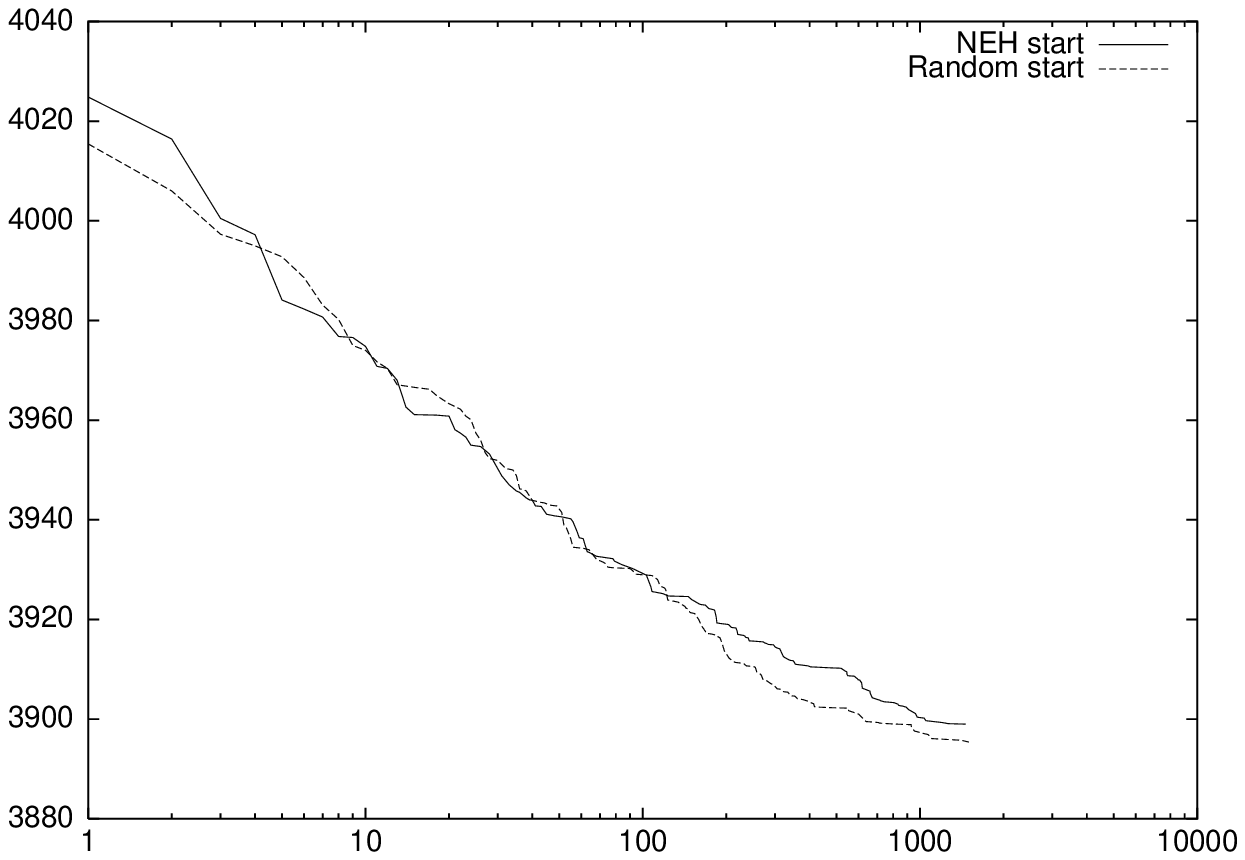}
	\includegraphics[width=0.5\textwidth]{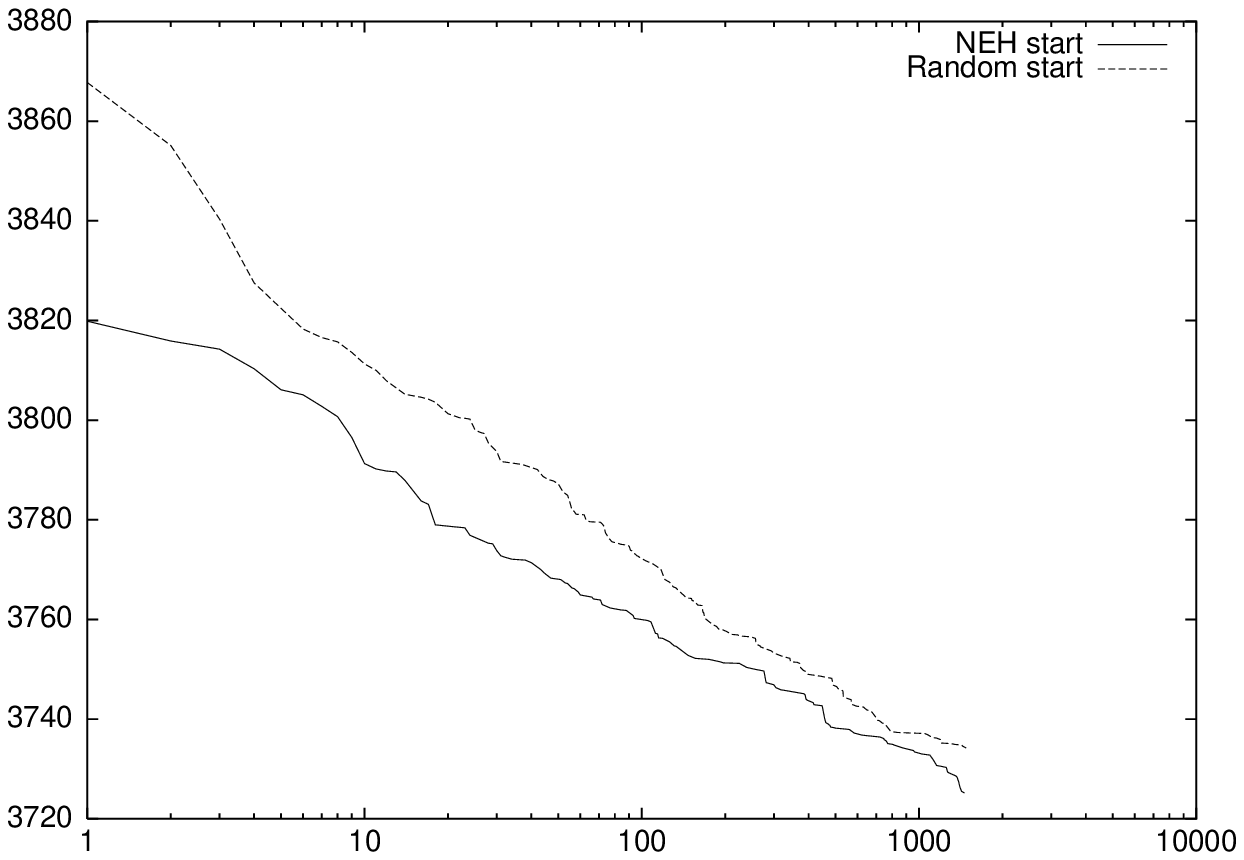}
\end{center}
\caption{The plots show the 
average solution quality (given on the $y$-axis) as a function of the
number of iterations (given on the $x$-axis) for an \ails\ algorithm
applied to the FSP on instances \inst{ta051} and \inst{ta056}. }
\label{f:fdcp}
\end{figure}

Let us illustrate the tradeoffs 
between random and greedy initial solutions when using an \ails\
algorithm for the permutation \fsp\ (FSP)~\cite{Stu98:aida98-04}. That
\ails\ algorithm uses a straight-forward local search implementation,
random perturbations, and always applies \genps\ to the best solution
found so far. In Figure~\ref{f:fdcp} we show how the average solution
cost evolves with the number of iterations for two instances.  The
averages are for $10$ independent runs when starting from random
initial solutions or from initial solutions returned by the NEH
heuristic~\cite{NawEnsHam83}. (NEH is one of the best performing
constructive heuristics for the FSP.) For short runs, the curve for
the instance on the right shows that
the NEH initial solutions lead to better average solution quality than
the random initial solutions. But at longer times, the picture is not
so clear, sometimes random initial solutions lead to better results as
we see on the instance on the left. This kind of test was also
performed for \ails\ applied to the TSP~\cite{AppCooRoh99}.
Again it was observed that the initial
solution had a significant influence on quality for short to medium 
sized runs.

In general, there will not always be a clear-cut answer regarding the best
choice of an initial solution, but greedy initial solutions appear to
be recommendable when one needs
low-cost solutions quickly. 
For much longer runs, the initial solution seems to be less relevant,
so the user can choose the initial solution that is the easiest
to implement. If however one has an application
where the influence of the initial
solution does persist for long times, probably the ILS walk is
having difficulty in exploring $\mathcal{S}^*$ and so
other perturbations or acceptance criteria should be considered.

\subsection{Perturbation}
\label{Sec3:ss:ps}

The main drawback of local descent is that it gets trapped in local
optima that are significantly worse than the global optimum.  Much
like simulated annealing, ILS escapes from local optima by applying
perturbations to the current local minimum. We will refer to the {\em
strength\/} of a perturbation as the number of solution components
which are modified. For instance for the TSP, it is the number of
edges that are changed in the tour, while in the \fsp, it is the
number of jobs which are moved in the perturbation. Generally, the
local search should not be able to undo the perturbation, otherwise
one will fall back into the local optimum just visited. Surprisingly
often, a {\it random} move in a neighborhood of higher order than the
one used by the local search algorithm can achieve this and will lead to a
satisfactory algorithm. Still better results can be obtained if the
perturbations take into account properties of the problem and are well
matched to the local search algorithm.

\begin{figure}
  \begin{center} 
	\includegraphics[width=0.5\textwidth]{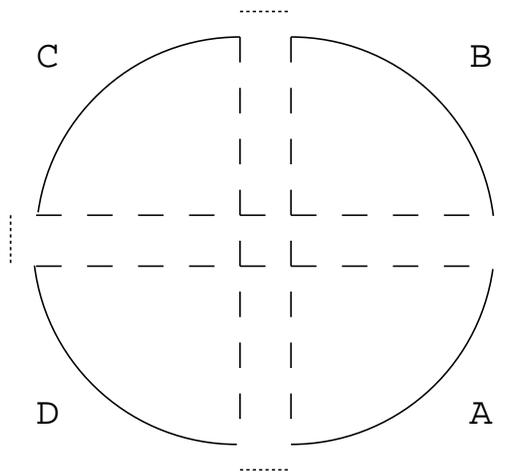}
  \end{center}
\caption{ Schematic representation of the double-bridge 
move. The four dotted edges are removed and the remaining parts
A, B, C, D are reconnected by the dashed edges.}
\label{f:db}
\end{figure}

By how much should the perturbation change the current solution? If
the perturbation is too strong, ILS may behave like a random restart,
so better solutions will only be found with a very low probability. On
the other hand, if the perturbation is too small, the local search
will often fall back into the local optimum just visited and the
diversification of the search space will be very limited. An example
of a simple but effective perturbation for the TSP is the {\em
double-bridge move\/}. This perturbation cuts four edges (and is thus
of ``strength'' $4$) and introduces four new ones as shown in
Figure~\ref{f:db}. Notice that each bridge is a $2$-change, but
neither of the $2$-changes individually keeps the tour connected.
Nearly all ILS studies of the TSP have incorporated this kind of
perturbation, and it has been found to be effective for all instance
sizes.  This is almost certainly because it changes the topology of
the tour and can operate on quadruples of very distant cities, whereas
local search always modifies the tour among nearby cities. (One could
imagine more powerful local searches which would include such
double-bridge changes, but the computational cost would be far greater
than for the local search methods used today.) In effect, the
double-bridge perturbation cannot be undone easily, neither by simple
local search algorithms such as \twoopt\ or \threeopt, nor by
Lin-Kernighan~\cite{LinKer73} which is currently the champion local
search algorithm for the TSP. Furthermore, this perturbation does not
increase much the tour length, so even if the current solution is very
good, one is almost sure the next one will be good, too.  These two
properties of the perturbation -- its small strength and its
fundamentally different nature from the changes used in local search
-- make the TSP the perfect application for \ils. But for other
problems, finding an effective perturbation may be more difficult.

We will now consider optimizing the perturbation
assuming the other modules to be fixed. In
problems like the TSP, one can hope to 
have a satisfactory ILS when using perturbations of fixed
size (independent of the instance size). On the contrary, for
more difficult problems, fixed-strength perturbations
may lead to poor performance.
Of course, the strength of the
perturbations used is not the whole story; their nature
is almost always very important and will also be discussed.
Finally we will close by pointing out that the perturbation strength
has an effect on the speed of the local search: weak
perturbations usually lead to faster execution of \genls. All these
different aspects need to be considered when optimizing 
this module.

\subsubsection{Perturbation strength}

For some problems, an appropriate perturbation strength is very small
and seems to be rather independent of the instance size. This is the
case for both the TSP and the FSP, and interestingly \ils\ 
for these problems is very competitive with today's 
best metaheuristic methods.
We can also consider other problems where instead one is driven
to large perturbation sizes.
Consider the example of an ILS algorithm for the \qap\
(QAP). We use an embedded \twoopt\ local search algorithm,
the perturbation is a random exchange of the location of $k$ items,
where $k$ is an adjustable parameter, and 
\genps\ always modifies the
best solution found so far. We applied this \ails\ algorithm to
QAPLIB instances\footnote{QAPLIB is accessible at
http://serv1.imm.dtu.dk/\~{}sk/qaplib/.} from four different classes of
QAP instances \cite{Tai95}; computational results are given in
Table~\ref{t:ilsqap}. A first observation is that the best
perturbation size is strongly dependent on the particular
instance. For two of the instances, the best performance was achieved
when as many as 75\% of the solution components were altered by the
perturbation. Additionally, for a too small perturbation strength, the
ILS performed worse than random restart (corresponding to the
perturbation strength $n$). However, the fact that random restart for
the QAP may perform---on average---better than a basic \ails\
algorithm is a bit misleading: in the next section we will show that
by simply modifying a bit the acceptance criterion, ILS becomes far
better than random restart. Thus one should keep in mind that the
optimization of an \ils\ may require more than the optimization of the
individual components.

\begin{table}[tb]
\caption{\small The first column is the name of the QAP instance;
the number gives its size $n$.  The successive columns are for
perturbation sizes $3$, $n/12$, $\cdots$, $n$.  A perturbation of size
$n$ corresponds to random restart. The table shows the mean solution
cost, averaged over $10$ independent runs for each instance. The
CPU-time for each trial is 30 sec.\ for \inst{kra30a}, 60 sec.\ for
\inst{tai60a} and \inst{sko64}, and 120 sec.\ for \inst{tai60b} on a 
Pentium III 500 MHz PC.}  
\label{t:ilsqap}
{\small
\begin{center}
\begin{tabular}{lrrrrrrrr}
instance & 3 & $n/12$ & $n/6$ & $n/4$ & $n/3$ & $n/2$ & $3n/4$ & $n$ \\\hline
\inst{kra30a} & 2.51 & 2.51 & 2.04 & 1.06 & 0.83 & 0.42 & 0.0  & 0.77\\
\inst{sko64}  & 0.65 & 1.04 & 0.50 & 0.37 & 0.29 & 0.29 & 0.82 & 0.93\\
\inst{tai60a} & 2.31 & 2.24 & 1.91 & 1.71 & 1.86 & 2.94 & 3.13 & 3.18\\
\inst{tai60b} & 2.44 & 0.97 & 0.67 & 0.96 & 0.82 & 0.50 & 0.14 & 0.43\\
\hline
\end{tabular}
\end{center}
}
\end{table}

\subsubsection{Adaptive perturbations}

The behavior of \ails\ for the QAP and also for other combinatorial
optimization problems \cite{HonKahMoo97,Stu98:aida98-04} shows that
there is no {\it \`a priori} single best size for the perturbation.
This motivates the possibility of modifying the perturbation strength
and adapting it {\it during} the run.

One possibility to do so is to exploit the search history. For the
development of such schemes, inspiration can be taken from what is
done in the context of \ts\ \cite{BatTec94:orsa,BatPro97:ea}. In
particular, Battiti and Protasi proposed~\cite{BatPro97:ea} a reactive
search algorithm for MAX-SAT which fits perfectly into the \ails\
framework. They perform a ``directed'' perturbation scheme which is
implemented by a tabu search algorithm and after each perturbation
they apply a standard local descent algorithm.

Another way of adapting the perturbation is to change
deterministically its strength during the search. One particular
example of such an approach is employed in the scheme called {\em
basic variable neighborhood search\/} (basic
VNS)~\cite{MlaHan97:cor,HanMla99:mic}; we refer to Section~\ref{s:r}
for some explanations on VNS. Other examples arise in the context of
\ts\ \cite{GloLag97}. In particular, ideas such as strategic
oscillations may be useful to derive more effective perturbations;
that is also the spirit of the reactive search algorithm previously
mentioned.

\subsubsection{More complex perturbation schemes}

Perturbations can be more complex than changes in a higher order
neighborhood. One rather general procedure to generate $s{'}$ from the
current $s^*$ is as follows. (1) Gently modify the definition of the
instance, e.g. via the parameters defining the various costs. (2) For
this modified instance, run \genls\ using $s^*$ as input; the output
is the perturbed solution $s{'}$. Interestingly, this is the method
proposed it the oldest ILS work we are aware of: in~\cite{Baxter81},
Baxter tested this approach with success on a location problem.  This
idea seems to have been rediscovered much later by Codenotti et al. in
the context of the TSP. Those authors~\cite{CodManMarRes96:informs}
first change slightly the city coordinates. Then they apply the local
search to $s^*$ using these perturbed city locations, obtaining the
new tour $s{'}$. Finally, running \genls\ on $s{'}$ using the {\it
unperturbed} city coordinates, they obtain the new candidate tour
$s^*{'}$.

Other sophisticated ways to generate good perturbations consist in
optimizing a sub-part of the problem. If this task is difficult for
the embedded heuristic, good results can follow. Such an approach was
proposed by Louren\c co~\cite{Ram95:ejor} in the context of the \jsp\
(JSP). Her perturbation schemes are based on defining one- or
two-machine sub-problems by fixing a number of variables in the
current solution and solving these sub-problems, either
heuristically~\cite{Lo98:Tech} or to optimality using for instance
Carlier's exact algorithm~\cite{Car82:ejor} or the early-late
algorithm \cite{Lo98:Tech}. These schemes work well because: (i) \ls\
is unable to undo the perturbations; (ii) after the perturbation, the
solutions tend to be very good and also have ``new'' parts that are
optimized.

\subsubsection{Speed}
\label{subsubsect_speed}

In the context of ``easy'' problems where ILS can work
very well with weak (fixed size) perturbations, there
is another reason why that metaheuristic can
perform much better than random restart: {\em Speed\/}. 
Indeed, \genls\ will usually execute
much faster on a solution obtained by applying a small perturbation to
a local optimum than on a random solution. As
a consequence, \ils\ can run many
more local searches than can random restart in the same CPU time.
As a qualitative example, consider
again Euclidean TSPs. ${\cal O}(n)$
local changes have to be applied by the local search
to reach a local optimum from a
random start, whereas empirically a nearly constant number is
necessary in ILS when using the $s{'}$ obtained with the
double-bridge perturbation. Hence, in a given amount of CPU time, ILS
can sample many more local optima than can random restart. This {\em
speed factor\/} can give ILS a considerable advantage over other restart
schemes.

Let us illustrate this speed factor quantitatively.
We compare for the TSP the number of local searches performed in a given
amount of CPU time by: (i) random restart; (ii) ILS using a 
double-bridge move; (iii) ILS using five
simultaneous double-bridge moves. (For both ILS implementations,
we used random starts and the routine \genac\ accepted only 
shorter tours.) For our numerical tests
we used a fast \threeopt\ implementation with standard speed-up
techniques. In particular, it used a fixed radius nearest neighbor
search within candidate lists of the 40 nearest neighbors for each
city and don't look bits
\cite{Ben92:orsa,JohMcG97,MarOttFel91:cs}. Initially, all don't look
bits are turned off (set to 0). If for a node no improving move can be
found, its don't look bit is turned on (set to 1) and the node is not
considered as a starting node for finding an improving move in the
next iteration. When an arc incident to a node is changed by a move,
the node's don't look bit is turned off again. 
In addition, when running ILS, after a
perturbation we only turn off the \dlbs\ of the 25 cities around each
of the four breakpoints in a current tour. All three algorithms were run for
120 seconds on a 266 MHz Pentium II processor on a set of
TSPLIB\footnote{TSPLIB is accessible at
www.iwr.uni-heidelberg.de/iwr/comopt/software/TSPLIB95.} instances
ranging from 100 up to 5915 cities. Results are given 
in Table~\ref{t:sf}. 
For small instances, we see that \ils\ ran between $2$ and $10$
times as many local searches as random restart. Furthermore, this advantage
of ILS grows fast with increasing
instance size: for the largest instance,
the first ILS algorithm ran approximately 260 times as many local searches as
random restart in our alloted time. 
Obviously, this speed advantage of ILS over random restart is strongly
dependent on the strength of the perturbation applied. The larger the
perturbation size, the more the solution is modified and generally the
longer the subsequent local search takes. This fact is intuitively
obvious and is confirmed in
Table~\ref{t:sf}.

\begin{table}[tb]
\label{t:sf}
\caption{\small The first column gives the name of the 
TSP instance which specifies its size. 
The next columns give the 
number of local searches performed when using: (i) random restart
($\#\mbox{\small LS}_{\tiny RR}$); (ii) ILS with a single double-bridge
perturbation ($\#\mbox{\small LS}_{\tiny 1-DB}$); (iii) 
ILS with a five
double-bridge perturbation
($\#\mbox{\small LS}_{\tiny 5-DB}$).
All algorithms were run 120 secs.\ on a Pentium 266 MHz PC.}  
{\small
\begin{center}
\begin{tabular}{lrrr}
instance & $\#\mbox{\small LS}_{\mbox{\tiny RR}}$ 
& $\#\mbox{\small LS}_{\mbox{\tiny 1-DB}}$ 
& $\#\mbox{\small LS}_{\mbox{\tiny 5-DB}}$ \\\hline
\inst{kroA100} & 17507 & 56186 & 34451 \\
\inst{d198}    & 7715  & 36849 & 16454 \\
\inst{lin318}  & 4271  & 25540 & 9430  \\
\inst{pcb442}  & 4394  & 40509 & 12880 \\
\inst{rat783}  & 1340  & 21937 & 4631  \\
\inst{pr1002}  & 910   & 17894 & 3345  \\
\inst{pcb1173} & 712   & 18999 & 3229  \\
\inst{d1291}   & 835   & 23842 & 4312  \\
\inst{fl1577}  & 742   & 22438 & 3915  \\
\inst{pr2392}  & 216   & 15324 & 1777  \\
\inst{pcb3038} & 121   & 13323 & 1232  \\
\inst{fl3795}  & 134   & 14478 & 1773  \\
\inst{rl5915}  & 34    & 8820  & 556   \\ 
\hline
\end{tabular}
\end{center}
}
\end{table}

In summary, the optimization of the perturbations depends
on many factors, and problem-specific characteristics play
a central role. Finally, it is important to keep in mind 
that the perturbations also interact with the other
components of \ails. We will discuss these interactions
in Section \ref{Sec3:ss:s}

\subsection{Acceptance criterion}
\label{Sec3:ss:ac}

ILS does a randomized walk in $\mathcal{S}^*$, the space of the local
minima. The perturbation mechanism together with the local search
defines the possible transitions between a current solution $s^*$ in
$\mathcal{S}^*$ to a ``neighboring'' solution $s^{*}{'}$ also in
$\mathcal{S}^*$. The procedure \genac\ then determines whether
$s^{*}{'}$ is accepted or not as the new current solution. \genac\ has
a strong influence on the nature and effectiveness of the walk in
$\mathcal{S}^*$.  Roughly, it can be used to control the balance
between intensification and 
diversification of that search. A simple way to illustrate
this is to consider a Markovian acceptance criterion. A very strong
intensification is achieved if only better solutions are accepted. We
call this acceptance criterion \ac{Better} and it is defined for
minimization problems as:

\begin{equation}
\label{e:better}
{\mbox{\ac{Better}}(s^*,s^*{'},\mbox{\em history\/})=}\left\{\matrix{
        \hfill s^*{'} \hskip 2pc  & \mbox{if}\ \mathcal{C}(s^*{'}) < \mathcal{C}(s^*) &\cr\cr
        \hfill s^* \hskip 2pc  & \mbox{otherwise}\hfill &
        }\right.
\end{equation}

At the opposite extreme is the random walk acceptance criterion (denoted by
\ac{RW}) which always applies the perturbation to the most recently
visited local optimum, irrespective of its cost: 

\begin{equation}
\label{e:rw}
{\mbox{\ac{RW}}}(s^*,s^*{'},\mbox{\em history\/})= s^*{'}
\end{equation}
This criterion clearly favors diversification over intensification.

Many intermediate choices between these two extreme cases are
possible. In one of the first ILS algorithms, the large-step Markov
chains algorithm proposed by Martin, Otto, and
Felten~\cite{MarOttFel91:cs,MarOttFel92:orl}, a simulated annealing
type acceptance criterion was applied. We call it
$\mbox{\ac{LSMC}}(s^*,s^{*}{'},\mbox{\em history\/})$.  In particular,
$s^*{'}$ is always accepted if it is better than $s^*$. Otherwise, if
$s^{*}{'}$ is worse than $s^*$, $s^{*}{'}$ is accepted with
probability $\exp\{(\mathcal{C}(s^*) -
\mathcal{C}(s^{*}{'})) / T\}$ where $T$ is a parameter called
temperature and it is usually lowered during the run as in
\sa. Note that \ac{LSMC} approaches the \ac{RW}
acceptance criterion if $T$ is very high, while at very low
temperatures \ac{LSMC} is similar to the \ac{Better} acceptance
criterion. An interesting possibility for \ac{LSMC} is to allow
non-monotonic temperature schedules as proposed in
\cite{HuKahTsa95} for simulated annealing or in tabu
thresholding \cite{Glo95}. This can be most effective if it is
done using memory: when further intensification no longer
seems useful, increase the temperature to do diversification for 
a limited time, then resume intensification. Of course, just as
in tabu search, it is desirable to do this in an automatic
and self-regulating manner~\cite{GloLag97}.

A limiting case of using memory in the acceptance criteria is to
completely restart the ILS algorithm when the intensification seems to
have become ineffective. (Of course this is a rather extreme way to
switch from intensification to diversification).  For instance one can
restart the ILS algorithm from a new initial solution if no improved
solution has been found for a given number of iterations. The restart
of the algorithm can easily be modeled by the acceptance criterion
called $\mbox{\ac{Restart}}(s^*,s^{*}{'},\mbox{\em history\/})$. Let
$i_{\mbox{\footnotesize\em last}}$ be the last iteration in which a
better solution has been found and $i$ be the iteration counter. Then
$\mbox{\ac{Restart}}(s^*,s^{*}{'},\mbox{\em history\/})$ is defined as

\begin{equation}
\label{c:ils:e:restart}
{\mbox{\ac{Restart}}(s^*,s^{*}{'},\mbox{\em history\/})=}\left\{\matrix{
        \hfill s^{*}{'} \hskip 2pc & \mbox{if}\ \mathcal{C}(s^{*}{'}) < \mathcal{C}(s^*)\hfill &\cr\cr
        \hfill s \hskip 2pc & \mbox{if}\ \mathcal{C}(s^{*}{'}) \geq \mathcal{C}(s^*) \mbox{ and } i -
        i_{\mbox{\footnotesize\em last}} > i_r &\cr\cr
        \hfill s^* \hskip 2pc  & \mbox{otherwise}.\hfill &
        }\right.
\end{equation}

\noindent where $i_r$ is a parameter that indicates that the
algorithm should be restarted if no improved solution was found for
$i_r$ iterations. Typically, $s$ can be generated in different
ways. The simplest strategy is to generate a new solution randomly or
by a greedy randomized heuristic. Clearly many other ways to
incorporate memory may and should be considered, the overall
efficiency of ILS being quite sensitive to the acceptance criterion
applied. We now illustrate this with two examples.

\subsubsection{Example 1: TSP}
Let us consider the effect of the two acceptance criteria
\ac{RW} and \ac{Better}. We performed our tests on the TSP as summarized in
Table \ref{t:sq}.  We give the average percentage excess over the
known optimal solutions when using $10$ independent runs on our set of
benchmark instances. In addition we also give this excess for the
random restart \threeopt\ algorithm.  First, we observe that both ILS
schemes lead to a significantly better average solution quality than
random restart using the same local search.  This is particularly true
for the largest instances, confirming again the claims given in
Section~\ref{sect_iterating} Second, given that one expects the good
solutions for the TSP to cluster (see Section~\ref{Sec3:ss:s}), a good
strategy should incorporate intensification. It is thus not surprising
to see that the \ac{Better} criterion leads to shorter tours than the
\ac{RW} criterion.

\begin{table}[tb]
\label{t:sq}
\caption{\small Influence of the acceptance criterion 
for various TSP instances. The first column gives the instance name
and its size. The next columns give the excess percentage length of
the tours obtained using: random restart (RR), \ils\ with \ac{RW}, and
\ils\ with \ac{Better}.  The data is averaged over $10$ independent
runs. All algorithms were run 120 secs.\ on a Pentium 266 MHz PC.}
{\small
\begin{center}
\begin{tabular}{lrrr}
instance & $\davg$(RR) & 
	$\davg$(\ac{RW}) & 
	$\davg$(\ac{Better}) \\\hline
\inst{kroA100} & 0.0  &  0.0 & 0.0 \\
\inst{d198}    & 0.003&  0.0 & 0.0 \\
\inst{lin318}  & 0.66 & 0.30 & 0.12 \\
\inst{pcb442}  & 0.83 & 0.42 & 0.11 \\
\inst{rat783}  & 2.46 & 1.37 & 0.12 \\
\inst{pr1002}  & 2.72 & 1.55 & 0.14 \\
\inst{pcb1173} & 3.12 & 1.63 & 0.40 \\
\inst{d1291}   & 2.21 & 0.59 & 0.28 \\
\inst{fl1577}  & 10.3 & 1.20 & 0.33 \\
\inst{pr2392}  & 4.38 & 2.29 & 0.54 \\
\inst{pcb3038} & 4.21 & 2.62 & 0.47 \\
\inst{fl3795}  & 38.8 & 1.87 & 0.58 \\
\inst{rl5915}  & 6.90 & 2.13 & 0.66 \\ 
\hline
\end{tabular}
\end{center}
}
\end{table}

The runs given in this example are rather short. For much longer runs,
the \ac{Better} strategy comes to a point where it no longer finds
improved tours and diversification should be considered again. Clearly it
will be possible to improve significantly the results by alternating
phases of intensification and diversification. 

\subsubsection{Example 2: QAP} Let us come back 
to ILS for the QAP discussed previously. For this problem we found
that the acceptance criterion \ac{Better} together with a (poor)
choice of the perturbation strength could result in worse performance
than random restart. In Table \ref{t:acqap} we give results for the
same ILS algorithm except that we now also consider the use of the
\ac{RW} and \ac{Restart} acceptance criteria. We see 
that the ILS algorithm using these modified acceptance criteria are
much better than random restart, the only exception being
\ac{RW} with a small perturbation strength on \inst{tai60b}. 

\begin{table}[tb]
\label{t:acqap}
\caption{\small Further tests on the QAP benchmark
problems using the same perturbations and CPU times as before; given is the mean solution
cost, averaged over $10$ independent runs for each instance. Here we
consider three different choices for the acceptance
criterion. Clearly, the inclusion of diversification significantly lowers
the mean cost found.}  {\small
\begin{center}
\begin{tabular}{llrrrrrrrr}
instance & acceptance & 3 & $n/12$ & $n/6$ & $n/4$ & $n/3$ & $n/2$ & $3n/4$ & $n$ \\\hline
\inst{kra30a} & \ac{Better} & 2.51 & 2.51 & 2.04 & 1.06 & 0.83 & 0.42 & 0.0  & 0.77\\
\inst{kra30a} & \ac{RW}     & 0.0  & 0.0  & 0.0  & 0.0  & 0.0  & 0.02 & 0.47 & 0.77\\
\inst{kra30a} & \ac{Restart}& 0.0  & 0.0  & 0.0  & 0.0  & 0.0  & 0.0  & 0.0  & 0.77\\
\hline
\inst{sko64}  & \ac{Better} & 0.65 & 1.04 & 0.50 & 0.37 & 0.29 & 0.29 & 0.82 & 0.93\\
\inst{sko64}  & \ac{RW}     & 0.11 & 0.14 & 0.17 & 0.24 & 0.44 & 0.62 & 0.88 & 0.93\\
\inst{sko64}  & \ac{Restart}& 0.37 & 0.31 & 0.14 & 0.14 & 0.15 & 0.41 & 0.79 & 0.93\\
\hline
\inst{tai60a} & \ac{Better} & 2.31 & 2.24 & 1.91 & 1.71 & 1.86 & 2.94 & 3.13 & 3.18\\
\inst{tai60a} & \ac{RW}     & 1.36 & 1.44 & 2.08 & 2.63 & 2.81 & 3.02 & 3.14 & 3.18\\
\inst{tai60a} & \ac{Restart}& 1.83 & 1.74 & 1.45 & 1.73 & 2.29 & 3.01 & 3.10 & 3.18\\
\hline
\inst{tai60b} & \ac{Better} & 2.44 & 0.97 & 0.67 & 0.96 & 0.82 & 0.50 & 0.14 & 0.43\\
\inst{tai60b} & \ac{RW}     & 0.79 & 0.80 & 0.52 & 0.21 & 0.08 & 0.14 & 0.28 & 0.43\\
\inst{tai60b} & \ac{Restart}& 0.08 & 0.08 & 0.005& 0.02 & 0.03 & 0.07 & 0.17 & 0.43\\
\hline
\end{tabular}
\end{center}
}
\end{table}

This example shows that there are strong inter-dependences between the
perturbation strength and the acceptance criterion. Rarely is this
inter-dependence completely understood. But, as a general rule of
thumb, when it is necessary to allow for diversification, we believe
it is best to do so by accepting numerous small perturbations rather
than by accepting one large perturbation.

Most of the acceptance criteria applied so far in ILS algorithms are
either fully Markovian or make use of the search history in a very
limited way. We expect that there will be many more ILS applications
in the future making strong use of the search history; in
particular, alternating between intensification and diversification is
likely to be an essential feature in these applications.

\subsection{Local search}
\label{Sec3:ss:ls}

So far we have treated the local search algorithm as a black box which
is called many times by \ails. Since the behavior and performance of
the over-all \ails\ algorithm is quite sensitive to the choice of the
embedded heuristic, one should optimize this choice whenever possible.
In practice, there may be many quite different algorithms that can be
used for the embedded heuristic. (As mentioned at the beginning of the
chapter, the heuristic need not even be a local search.) One might think
that the better the \ls, the better the corresponding
\ails. Often this is true. For instance in the context of the
TSP, Lin-Kernighan~\cite{LinKer73} is a better \ls\ than $3$-opt which
itself is better than $2$-opt \cite{JohMcG97}. Using a fixed type of
perturbation such as the double-bridge move, one finds that iterated
Lin-Kernighan gives better solutions than iterated $3$-opt which
itself gives better solutions than iterated $2$-opt
\cite{JohMcG97,StuHoo00:iridia}. But if we assume that the total
computation time is fixed, it might be better to apply more frequently
a faster but less effective local search algorithm than a slower and
more powerful one. Clearly which choice is best depends on just how
much more time is needed to run the better heuristic. If the speed
difference is not large, for instance if it is independent of the
instance size, then it usually worth using the better heuristic. This
is the most frequent case; for instance in the TSP, 3-opt is a bit
slower than 2-opt, but the improvement in quality of the tours are
well worth the extra CPU time, be-it using random restart or \ils. The
same comparison applies to using L-K rather than 3-opt. However, there
are other cases where the increase in CPU time is so large compared to
the improvement in solution quality that it is best not to use the
``better'' local search. For example, again in the context of the TSP,
it is known that $4$-opt gives slightly better solutions than $3$-opt, but
in standard implementations it is $O(n)$ times slower ($n$
being the number of cities). It is then better not to use $4$-opt as
the \ls\ embedded in ILS.\footnote{But see ref.~\cite{Glo96} for a way
to implement $4$-opt much faster.}

There are also other aspects that should be considered when selecting
a \ls. Clearly, there is not much point in having an excellent \ls\ if
it will systematically undo the perturbation; however this issue is
one of globally optimizing \ils, so it will be postponed till the next
sub-section. Another important aspect is whether one can really get
the speed-ups that were mentioned in sub-section~\ref{Sec3:ss:ps}.
There we saw that a standard trick for \genls\ was to introduce
\dlbs. These give a large gain in speed if the bits can be reset also
after the application of the perturbation. This requires that the
developper be able to access the source code of \genls. A state of the
art ILS will take advantage of all possible speed-up tricks, and thus
the \genls\ most likely will not be a true black box.

Finally, there may be some advantages in allowing \genls\ to sometimes
generate worse solutions. For instance, if we replace the \ls\
heuristic by tabu search or short simulated annealing runs, the corresponding
\ails\ may perform better. This seems most promising when standard \ls\
methods perform poorly. Such is indeed the case in the job-shop
scheduling problem: the use of tabu search as the embedded heuristic
gives rise to a very effective iterated local
search~\cite{RamZwi96:mic}.

\subsection{Global optimization of ILS}
\label{Sec3:ss:s}

So far, we have considered representative issues arising when
optimizing separately each of the four components of an \ils. In
particular, when illustrating various important characteristics of one
component, we kept the other components fixed. But clearly the
optimization of one component depends on the choices made for the
others; as an example, we made it clear that a good perturbation must
have the property that it cannot be easily undone by the \ls. Thus, at
least in principle, one should tackle the {\it global} optimization of
an ILS.  Since at present there is no theory for analyzing a
metaheuristic such as \ils, we content ourselves here with just giving
a rough idea of how such a global optimization can be approached in
practice.

If we reconsider the sub-section on the effect of the initial
solution, we see that \genis\ is to a large extent irrelevant when the
ILS performs well and rapidly looses the memory of its starting point.
Hereafter we assume that this is the case; then the optimization of
\genis\ can be ignored and we are left with the joint optimization of
the three other components. Clearly the best choice of \genps\ depends
on the choice of \genls\ while the best choice of \genac\ depends on
the choices of \genls\ and \genps\/.  In practice, we can approximate
this global optimization problem by successively optimizing each
component, assuming the others are fixed until no improvements are
found for any of the components. Thus the only difference with what
has been presented in the previous sub-sections is that the
optimization has to be iterative. This does not guarantee global
optimization of the ILS, but it should lead to an adequate
optimization of the overall algorithm.

Given these approximations, we should make more precise what in fact
we are to optimize.  For most users, it will be
the mean (over starting solutions) of the best cost found during a run
of a given length. Then the ``best'' choice for the different
components is a well posed problem, though it is intractable without
further restrictions. Furthermore, in general, the detailed instance
that will be considered by the user is not known ahead of time, so it
is important that the resulting ILS algorithm be robust. Thus it is
preferable not to optimize it to the point where it is sensitive to
the details of the instance. This robustness seems to be achieved in
practice: researchers implement versions of \ils\ with a reasonable level of
global optimization, and then test with some degree of success the
performance on standard benchmarks.

\subsubsection*{Search space characteristics.}
\label{sss:ssc}

At the risk of repeating ourselves, let us highlight the main
dependencies of the components:
\begin{enumerate}
\item the perturbation should not be easily undone by the \ls;
if the \ls\ has obvious short-comings, a good perturbation
should compensate for them.
\item the combination \genps--\genac\ determines the relative
balance of intensification and diversification; large perturbations 
are only useful if they can be accepted, which occurs only
if the acceptance criterion is not too biased towards
better solutions.
\end{enumerate} 
As a general guideline, \genls\ should be as powerful as possible as
long as it is not too costly in CPU time. Given such a choice, then
find a well adapted perturbation following the discussion in
Section~\ref{Sec3:ss:ps}; to the extent possible, take advantage of
the structure of the problem. Finally, set the
\genac\ routine so that ${\mathcal S}^*$ is sampled
adequately. With this point of view, the overall optimization of the
ILS is nearly a bottom-up process, but with iteration. Perhaps the
core issue is what to put into \genps: can one restrict the
perturbations to be weak? From a theoretical point of view, the answer
to this question depends on whether the best solutions ``cluster'' in
${\cal S}^*$. In some problems (and the TSP is one of them), there is
a strong correlation between the cost of a solution and its
``distance'' to the optimum: in effect, the best solutions cluster
together, i.e., have many similar components.  This has been referred
to in many different ways: ``Massif Central''
phenomenon~\cite{FonRobPreTal99:mic}, principle of proximate
optimality~\cite{GloLag97}, and replica
symmetry~\cite{MezardParisi87b}. If the problem under consideration
has this property, it is not unreasonable to hope to find the true
optimum using a biased sampling of ${\cal S}^*$. In particular, it is
clear that is useful to use intensification to improve the probability
of hitting the global optimum.

There are, however, other types of problems where the clustering is
incomplete, i.e., where very distant solutions can be nearly as good
as the optimum.  Examples of combinatorial optimization problems in
this category are QAP, graph bi-section, and MAX-SAT. When the space
of solutions has this property, new strategies have to be
used. Clearly, it is still necessary to use intensification to get the
best solution in one's current neighborhood, but generally this will
not lead to the optimum.  After an intensification phase, one must go
explore other regions of ${\cal S}^*$. This can be attempted by using
``large'' perturbations whose strength grows with the instance. Other
possibilities are to restart the algorithm from scratch and repeat
another intensification phase or by oscillating the acceptance
criterion between intensification and diversification
phases. Additional ideas on the tradeoffs between intensification and
diversification are well discussed in the context of \ts\ (see, for
example, \cite{GloLag97}). Clearly, the balance intensification --
diversification is very important and is a challenging problem.

\section{Selected applications of ILS}
\label{s:rw}

ILS algorithms have been applied successfully to a variety of
combinatorial optimization problems. In some cases, these algorithms
achieve extremely high performance and even constitute the current
state-of-the-art metaheuristics, while in other cases the ILS approach
is merely competitive with other metaheuristics. In this section, we
cover some of the most studied problems, with a stress on the \tsp\
and scheduling problems.

\subsection{ILS for the TSP}

The TSP is probably the best-known combinatorial optimization problem.
{\it De facto}, it is a standard test-bed for the development of new
algorithmic ideas: a good performance on the TSP is taken as evidence
of the value of such ideas. Like for many metaheuristic algorithms,
some of the first ILS algorithms were introduced and tested on the
TSP, the oldest case of this being due to
Baum~\cite{Baum_86a,Baum_86b}. He coined his method {\em iterated
descent\/}; his tests used 2-opt as the embedded heuristic, random
3-changes as the perturbations, and imposed the tour length to
decrease (thus the name of the method).  His results were not
impressive, in part because he considered the non-Euclidean TSP, which
is substantially more difficult in practice than the Euclidean TSP. A
major improvement in the performance of ILS algorithms came from the
{\em large-step Markov chain\/} (LSMC) algorithm proposed by Martin,
Otto, and Felten \cite{MarOttFel91:cs}. They used a simulated
annealing like acceptance criterion (\ac{LSMC}) from which the
algorithm's name is derived and considered both the application of
\threeopt\ local search and the Lin-Kernighan heuristic (\lk ) which is
the best performing local search algorithm for the TSP. But probably
the key ingredient of their work is the introduction of the
double-bridge move for the perturbation. This choice made the approach
very powerful for the Euclidean TSP, and that encouraged much more
work along these lines.  In particular, Johnson~\cite{Joh90,JohMcG97}
coined the term ``iterated Lin-Kernighan'' (ILK) for his
implementation of ILS using the Lin-Kernighan as the local search.
The main differences with the LSMC implementation are: (i)
double-bridge moves are random rather than biased; (ii) the costs
are improving (only better tours are accepted, corresponding to the
choice \ac{Better} in our notation).  Since these initial studies,
other ILS variants have been proposed, and Johnson and
McGeoch~\cite{JohMcG97} give a summary of the situation as of 1997.

Currently the highest performance \ails\ for the TSP is the chained
\lk\ code by Applegate, Bixby, Chvatal, and Cook which is available as
a part of the Concorde software package at
www.keck.caam.rice.edu/\-concorde.\-html. These authors have provided
very detailed descriptions of their implementation, and so we refer
the reader to their latest article~\cite{AppBixChvCoo00} for
details. Furthermore, Applegate, Cook, and Rohe~\cite{AppCooRoh99}
performed thorough experimental tests of this code by considering the
effect of the different modules: (i) initial tour; (ii) implementation
choices of the \lk\ heuristic; (iii) types of perturbations.  Their
tests were performed on very large instances with up to 25 million
cities. For the double-bridge move, they considered the effect of
forcing the edges involved to be ``short'', and investigated 
the random double-bridge moves as well.
Their conclusion is that the best performance is obtained when
the double-bridge moves are biased towards short edge
lengths. However, the strength of the bias towards short edges should
be adapted to the available computation time: the shorter the
computation time, the shorter the edges should be. In their tests on
the influence of the initial tour, they concluded that the worst
performance is obtained with random initial tours or those returned by
the nearest neighbor heuristic, while best results were obtained with
the Christofides algorithm \cite{Chr76}, the greedy heuristic
\cite{Ben92:orsa} or the Quick-Boruvka heuristic proposed in that article.
With long runs of their algorithm on TSPLIB instances with more than 10.000
cities they obtained an impressive performance, always obtaining
solutions that have less than 0.3\% excess length over the lower bounds for
these instances. For the largest instance considered, a 25 million city
instance, they reached a solution of only 0.3\% over the estimated optimum.

Apart from these works, two new ILS algorithms for the TSP have been
proposed since the review article of Johnson and McGeoch. The first
algorithm is due to St\"utzle~\cite{Stu98:phd,StuHoo00:iridia}; he
examined the run-time behavior of ILS algorithms for the TSP and
concluded that ILS algorithms with the \ac{Better} acceptance
criterion show a type of stagnation behavior for long
run-times~\cite{Stu98:phd} as expected when performing a strong
intensification search.  To avoid such stagnation, restarts and a
particular acceptance criterion to diversify the search were proposed.
The goal of this latter strategy is to force the search to continue
from a position that is beyond a certain minimal distance from the
current position.  This idea is implemented as follows. Let $s_c$ be
the solution from which to escape; $s_c$ is typically chosen as
$s^*_{\mbox{\footnotesize\em best}}$, the best solution found in the
recent search. Let $d(s,s')$ be the distance between two tours $s$
and $s'$, that is the number of edges in which they differ. Then 
the following steps
are repeated until a solution beyond a minimal distance $d_{\min}$
from $s_c$ is obtained:

\begin{itemize}
\item[$(1)$] Generate $p$ copies of $s_c$.
\item[$(2)$] To each of the $p$ solutions apply \genps\ followed by \genls .
\item[$(3)$] Choose the best $q$ solutions, $1 < q \leq p$, as candidate
  solutions.
\item[$(4)$] Let $s^*$ be the candidate solution with maximal distance to
  $s_c$.  If $d(s^*,s_c) \leq d_{\min}$ then repeat at $(2)$; otherwise
  return $s^*$.
\end{itemize}

The purpose of step 3 is to choose good quality solutions, while step
4 guarantees that the point from which the search will be continued is
sufficiently different (far) from $s_c$. The attempts are continued
until a new solution is accepted, but one gives up after some maximum
number of iterations. Computational results for this way of
going back and forth
between intensification and diversification show that the
method is very effective, even when using only a \threeopt\ local
search~\cite{StuHoo00:iridia,StuGruLinRut00:ppsn}.

The second \ails\ developed for the TSP since 1997 is that of
Katayama and Narisha~\cite{KatNar99:gecco}. They introduce a new
perturbation mechanism which they called a {\it genetic
transformation}. The genetic transformation mechanism uses two tours,
one of which is the best found so far $s^*_{\mbox{\footnotesize\em
best}}$, while the second solution $s'$ is a tour found earlier in the
search. First a random 4-opt move is performed on
$s^*_{\mbox{\footnotesize\em best}}$, resulting in $s^{*}{'}$. Then
the subtours that are shared among $s^{*}{'}$ and $s'$ are enumerated.
The resulting parts are then reconnected with a greedy
algorithm. Computational experiments with an iterated \lk\ algorithm
using the genetic transformation method instead of the standard
double-bridge move have shown that the approach is very effective;
further studies should be forthcoming.

\subsection{ILS for scheduling problems}

ILS has also been applied successfully to scheduling problems. Here we
summarize the different uses of ILS for tackling these types of
systems, ranging from single machine to complex multi-machine
scheduling.

\subsubsection{Single Machine Total Weighted Tardiness Problem (SMTWTP)} 
Congram, Potts and van de Velde~\cite{ConPotVel98} have presented an
ILS algorithm for the SMTWTP based on a dynasearch local
search. Dynasearch uses dynamic programming to find a best move which
is composed of a set of independent interchange moves; each such move
exchanges the jobs at positions $i$ and $j$, $j\neq i$. Two
interchange moves are independent if they do not overlap, that is if
for two moves involving positions $i,j$ and $k,l$ we have
$\min\{i,j\} \geq \max \{k,l\}$ or vice versa. This neighborhood is of
exponential size but dynasearch explores this neighborhood in
polynomial time.

The perturbation consists of a series of random interchange
moves. They also exploit a well-known property of the SMTWTP: there
exists an optimal solution in which non-late jobs are sequenced in
non-decreasing order of the due dates. This property is used in two
ways: to diversify the search in the perturbation step and to reduce
the computation time of the dynasearch. In the acceptance criterion,
Congram et al.\  introduce a {\em backtrack step\/}: after $\beta$
iterations in which every new local optimum is accepted, the algorithm
restarts with the best solution found so far. In our notation, the
backtrack step is a particular choice for the history dependence incorporated
into \genac.

Congram et al.\ used several different embedded \genls, all associated
with the interchange neighborhood. These heuristics were: (i)
dynasearch; (ii) a \ls\ based on first-improvement descent; (iii) a
\ls\ based on best-improvement descent.  Then they performed tests to
evaluate these algorithms using random restart and compared them to
using \ils. While random restart dynasearch performed only slightly
better than the two simpler descent methods, the \ails\ with
dynasearch significantly outperformed the other two iterated descent
algorithms, which in turn were far superior to the random restart
versions.  The authors also show that the iterated dynasearch
algorithm significantly improves over the previously best known
algorithm, a tabu search presented in~\cite{CraPoWa98:IJC}.

\subsubsection{Single and parallel machine scheduling}
Brucker, Hurink, and Werner~\cite{BruHurWer96:dam,BruHurWer97:dam}
apply the principles of ILS to a number of one-machine and
parallel-machine scheduling problems. They introduce a local search
method which is based on two types of neighborhoods. At each step one
goes from one feasible solution to a neighboring one with respect to
the secondary neighborhood. The main difference with standard local
search methods is that this secondary neighborhood is defined on the
set of locally optimal solutions with respect to the first
neighborhood. Thus in fact this is an \ails\ with two nested
neighborhoods; searching in their primary neighborhood corresponds to
our local search phase; searching in their secondary neighborhood is
like our perturbation phase. The authors also note that the second
neighborhood is problem specific; this is what arises in \ails\ where
the perturbation should be adapted to the problem. The search at a
higher level reduces the search space and at the same time leads to
better results.

\subsubsection{Flow shop scheduling}

St\"utzle~\cite{Stu98:aida98-04} applied ILS to the \fsp\ (\afsp
). The algorithm is based on a straightforward first-improvement local
search using the insert neighborhood, where a job at position $i$ is
removed and inserted at position $j\neq i$. The initial schedule is
constructed by the NEH heuristic~\cite{NawEnsHam83} while the
perturbation is generated by composing moves of two different kinds:
swaps which exchange the positions of two adjacent jobs, and
interchange moves which have no constraint on adjacency.
Experimentally, it was found that perturbations with just a few swap
and interchange moves were sufficient to obtain very good results. The
article also compares different acceptance criteria;
\ac{ConstTemp}, which is the same as the \ac{LSMC} acceptance criterion except
that it uses a constant temperature $T_c$, was found to be superior to
\ac{Better}. The computational results show that despite the
simplicity of the approach, the quality of the solutions obtained is
comparable to that of the best performing local search algorithms for
the FSP; we refer to~\cite{Stu98:aida98-04} for a more detailed
discussion.

ILS has also been used to solve a flow-shop problem with several
stages in series.  Yang, Kreipl and Pinedo~\cite{YaKrPi00:js}
presented such a method; at each stage, instead of a single machine,
there is a group of identical parallel machines. Their metaheuristic
has two phases that are repeated iteratively. In the first phase, the
operations are assigned to the machines and an initial sequence is
constructed. The second phase uses an ILS to find better schedules for
each machine at each stage by modifying the sequence of each
machine. (This part is very similar in spirit to the approach of
Kreipl for the minimum total weighted tardiness job-shop
problem~\cite{Kreipl00:js} that is presented below.)  Yang,
Kreipl and Pinedo also proposed a ``hybrid'' metaheuristic: they first
apply a decomposition procedure that solves a series of single stage
sub-problems; then they follow this by their ILS. The process is
repeated until a satisfactory solution is obtained.

\subsubsection{Job shop scheduling}

Louren\c{c}o~\cite{Ram95:ejor} and Louren\c{c}o and
Zwijnenburg~\cite{RamZwi96:mic} used ILS to tackle the \jsp\ (\ajsp
). They performed extensive computational tests, comparing different
ways to generate initial solutions, various local search algorithms,
different perturbations, and three acceptance criteria. While they
found that the initial solution had only a very limited influence, the
other components turned out to be 
very important. Perhaps the heart of their work
is the way they perform the perturbations. They consider relaxations
of the problem at hand corresponding to the optimization of just some
of the jobs. Then they use exact methods to solve these sub-problems,
generating the perturbation move. This has the great advantage that
much problem-specific knowledge is built into the perturbation.  Such
problem specific perturbations are
difficult to generate from local moves only. Now, for the local search,
three alternatives were considered: local descent, short simulated
annealing runs, and short tabu search runs. Best results were obtained
using the latter in the local search phase. Not surprisingly, ILS
performed better than random restart given the same amount of time,
for any choice of the embedded local search heuristic.

In more recent work on the job-shop scheduling problem, Balas and
Vazacopoulos~\cite{BaVa98:MS} presented a variable depth search
heuristic which they called guided local search (GLS). GLS is based on
the concept of neighborhood trees, proposed by the authors, where each
node corresponds to a solution and the child nodes are obtained by
performing an interchange on some critical arc. In their work, the
interchange move consists in reversing more than one arc and can be
seen as a particular kind of variable depth interchange. They
developed ILS algorithms by embedding GLS within the shifting
bottleneck (SB) procedure by replacing the reoptimization cycle of the
SB with a number of cycles of the GLS procedure. They call
this procedure SB-GLS1. Later, they also proposed a variant of this
method, SB-GLS2, which works as follows. After all machines have been
sequenced, they iteratively remove one machine and apply GLS to a
smaller instance defined by the remaining machines. Then again GLS is
applied on the initial instance containing {\em all\/} machines. This
procedure is an ILS where a perturbed solution is obtained by applying
a (variable depth) local search to just part of an instance. 
The authors perform a
computational comparison with other metaheuristics and conclude that
SB-GLS (1 and 2) are robust and efficient, and provide schedules of
high quality in a reasonable computing time.
In some
sense, both heuristics are similar to the one proposed by
Louren\c{c}o~\cite{Ram95:ejor}, the main differences being: (i)
Louren\c co's heuristic applies perturbations to complete schedules whereas
the SB-GLS heuristic starts by an empty (infeasible) schedule and
iteratively optimizes it machine by machine until all machines have been
scheduled, in a SB-style followed by a local search application;
(ii) the local search algorithms used differ.

Recently, Kreipl applied ILS to the total weighted tardiness job shop
scheduling problem (TWTJSP) \cite{Kreipl00:js}. The TWTJSP is closer
to real life problems than the classical JSP with makespan objective
because it takes into account release and due dates and also
it introduces weights
that indicate the importance of each job. Kreipl uses an ILS algorithm
with the \ac{RW} acceptance criterion. The algorithm starts with an
initial solution obtained by the shortest processing time
rule~\cite{Haupt89:ORS}. The local search 
consists in reversing critical arcs and arcs adjacent to these,
where a critical arc has to be an element of at least one
critical path (there may exist several critical paths). One original
aspect of this ILS is the perturbation step: Kreipl applies a few
steps of a simulated annealing type algorithm with the Metropolis
acceptance criterion~\cite{Metropolis53:JCP} but with a fixed
temperature. For this perturbation phase a smaller neighborhood
than the one used in the local search phase is taken: while in the
local search phase any critical arc can be reversed, during the
diversification phase only the critical arcs belonging to the critical
path having the job with highest impact on the objective function are
considered.\footnote{It should be noted that the perturbation phase
leads, in general, to an intermediate solution which is not locally
optimal.} The number of iterations performed in the perturbation
depends how good the incumbent solution is. In promising regions, only
a few steps are applied to stay near good solutions, otherwise, a
"large" perturbation is applied to permit the algorithm to escape from
a poor region. Computational results with the ILS algorithm on a set
of benchmark instances has shown a very promising performance compared
to an earlier shifting bottleneck heuristic~\cite{Singer97:IIESL}
proposed for the same problem.

\subsection{ILS for other problems}

\subsubsection{Graph bipartitioning} 
\ails\ algorithms have been proposed and tested
on a number of other problems, though not as thoroughly as the ones we
have discussed so far.  We consider first the graph bipartitioning
problem. Given a (weighted) graph and a bisection or partition of its
vertices into two sets $A$ and $B$ of equal size, call the cut of the
partition the sum of the weights of the edges connecting the two
parts. The graph partitioning problem is to find the partition with
the minimum cut.  Martin and Otto~\cite{MarOtt95:cpe,MarOtt96:aor}
introduced an \ails\ for this problem following their earlier work on
the TSP.  For the local search, they used the Kernighan-Lin variable
depth local search algorithm (\kl ) \cite{KerLin70} which is the
analog for this problem of the \lk\ algorithm. In effect,
\kl\ finds intelligently $m$ vertices of one set to be exchanged
with $m$ of the other. Then, when considering possible perturbations,
they noticed a particular weakness of the \kl\ local search: \kl\
frequently generates partitions with many ``islands'', i.e., the two
sets $A$ and $B$ are typically highly fragmented (disconnected). Thus
they introduced perturbations that exchanged vertices between these
islands rather than between the whole sets $A$ and $B$.  This works as
follows: choose at random one of the cut edges, i.e., an edge
connecting $A$ and $B$. This edge connects two ``seed'' vertices each
belonging to their island. Around each seed, iteratively grow a
connected cluster of vertices within each island. When a target
cluster size or a whole island size is reached, stop the growth. The
two clusters are then exchanged and this is the perturbation
move. Finally, for the acceptance criterion, Martin and Otto used the
\ac{Better} acceptance criterion. The overall algorithm significantly
improved over the embedded local search (random restart of \kl ); it
also improved over simulated annealing if the acceptance criterion was
optimized.

At the time of that work, simulated annealing was the state of the art
method for the graph bisection problem. Since then, there have been
many other metaheuristics~\cite{BatBer99,MerFre00:ec} developed for
this problem, so the performance that must be reached is much higher
now.  Furthermore, given that the graph bipartitioning problem has a
low cost-distance correlation~\cite{MerFre00:ec}, ILS has difficulty
in sampling all good low cost solutions. To overcome this, some form of
history dependence most certainly would have to be built into
the perturbation or the acceptance criterion.

\subsubsection{MAX-SAT}
Battiti and Protasi present an application of {\em reactive search\/}
to the MAX-SAT problem~\cite{BatPro97:ea}. Their algorithm consists of
two phases: a local search phase and a diversification (perturbation)
phase.  Because of this, their approach fits perfectly into the ILS
framework. Their perturbation is obtained by running a tabu search on
the current local minimum so as to guarantee that the modified
solution $s'$ is sufficiently different from the current solution
$s^*$. Their measure of difference is just the Hamming distance; the
minimum distance is set by the length of a tabu list that is adjusted
during the run of the algorithm. For the \genls, they use a standard
greedy descent local search appropriate 
for the MAX-SAT problem.  Depending on the
distance between $s^{*}{'}$ and $s^*$, the tabu list length for the
perturbation phase is dynamically adjusted. The next perturbation
phase is then started based on solution $s^{*}{'}$ --- corresponding to
the \ac{RW} acceptance criterion. This work illustrates very nicely
how one can adjust dynamically the perturbation strength in an \ails\
run.  We conjecture that similar schemes will prove useful to optimize
\ails\ algorithms in a nearly automatic way.

\subsubsection{Prize-collecting Steiner tree problem}
The last combinatorial optimization problem we discuss is the
prize-collecting Steiner tree problem on graphs.  Canudo, Resende and
Ribeiro \cite{Resende00:tech} presented several local search
strategies for this problem: iterative improvement, multi-start with
perturbations, path-relinking, variable neighborhood search, and a
algorithm based on the integration of all these. They showed that all
these strategies are effective in improving solutions; in fact in many
of their tests they found the optimal solution. One of their proposed
heuristics, local search with perturbations, is in fact an ILS. In
that approach, they first generated initial solutions by the
primal-dual algorithm of Goemans and Wiliamson (GW)~\cite{GoeWi96} but
where the cost function is slightly modified.  Canudo et al.\ proposed
two perturbation schemes: perturbation by eliminations and
perturbations by prize changes. In the first scheme, the perturbation
is done by resetting to zero the prizes of some persistent node which
appeared in the solution build by GW and remained at the end of local
search in the previous iteration.  In the second scheme, the
perturbation consists in introducing noise into the node
prize. This feature of always applying the perturbation to the last
solution obtained by the local search phase is clearly in our notation
the ILS-RW choice.

\subsection{Summary}

The examples we have chosen in this section stress several points that
have already been mentioned before.  First, the choice of the local
search algorithm is usually quite critical if one is to obtain peak
performance. In most of the applications, the best performing \ails\
algorithms apply much more sophisticated local search algorithms than
simple best- or first-improvement descent methods.  Second, the other
components of an ILS also need to be optimized if the state of the art
is to be achieved. This optimization should be global, and to succeed
should involve the use of problem-specific properties. Examples of
this last point were given for instance in the scheduling
applications: there the good perturbations were not simply random
moves, rather they involved re-optimizations of significant parts of
the instance (c.f. the job shop case).

The final picture we reach is one where (i) \ails\ is a versatile
metaheuristic which can easily be adapted to different combinatorial
optimization problems; (ii) sophisticated perturbation schemes and
search space diversification are the essential ingredients to achieve the
best possible ILS performance.

\section{Relation to other metaheuristics}
\label{s:r}

In this section we highlight the similarities and differences between
\ails\ and other well-known metaheuristics.  We shall distinguish
metaheuristics which are essentially variants of \ls\ and those which
generate solutions using a mechanism that is not necessarily based on
an explicit neighborhood structure. Among the first class which we
call {\it neighborhood based metaheuristics} are methods like \sa\
(\asa )~\cite{KirGelVec83,Cer85}, \ts\ (\ats
)~\cite{Glo89,Glo90,GloLag97} or guided local search
(GLS)~\cite{VouTsa95}. The second class comprises metaheuristics like
GRASP~\cite{FeoRes95}, ant colony optimization
(\aaco)~\cite{DorDic99:nio}, evolutionary algorithms (EA)
\cite{Bae96,MicFog99}, scatter search~\cite{Glo99:nio}, \vns\ (\avns
)~\cite{HanMla99:mic,MlaHan97:cor} and \ails . Some metaheuristics of
this second class, like EAs and ACO, do not necessarily make use of
local search algorithms; however a \ls\ can be embedded in them, in
which case the performance is usually enhanced \cite{Muh91,Stu98:phd}.
The other metaheuristics in this class explicitly use embedded local
search algorithms as an essential part of their structure. For
simplicity, we will assume in what follows that all the metaheuristics
of this second class do incorporate local search algorithms. In this
case, such metaheuristics generate iteratively input solutions that
are passed to a local search; they can thus be interpreted as
multi-start algorithms, using the most general meaning of that term.
This is why we call them here {\it multi-start based
metaheuristics\/}.

\subsection{Neighborhood based metaheuristics}

Neighborhood based metaheuristics are extensions of iterative
improvement algorithms and avoid getting stuck in locally optimal
solutions by allowing moves to worse solutions in one's neighborhood.
Different metaheuristics of this class differ mainly by their move
strategies. In the case of \sa, the neighborhood is sampled randomly
and worse solutions are accepted with a probability which depends on a
temperature parameter and the degree of deterioration incurred; better
neighboring solutions are usually accepted while much worse
neighboring solutions are accepted with a low probability. In the case
of (simple) \ts\ strategies, the neighborhood is explored in an
aggressive way and cycles are avoided by declaring attributes of visited
solutions as tabu. Finally, in the case of guided local search, the
evaluation function is dynamically modified by penalizing certain
solution components. This allows the search to escape from a solution
that is a local optimum of the original objective function.

Obviously, any of these neighborhood based metaheuristics can be used
as the \genls\ procedure in \ails. In general, however, those
metaheuristics do not halt, so it is necessary to limit their run time if
they are to be embedded in \ails. One particular advantage of
combining neighborhood based metaheuristics with \ails\ is that they
often obtain much better solutions than iterative descent
algorithms. But this advantage usually comes at the cost of larger
computation times. Since these metaheuristics allow one to obtain
better solutions at the expense of greater computation times, we are
confronted with the following optimization problem when using them
within an \ails:
\footnote{This question is not specific to \ails; it arises for all
multi-start type metaheuristics.} ``For how long should one run the
embedded search 
in order to achieve the best tradeoff between computation
time and solution quality?'' This is very analogous to the question of
whether it is best to have a fast but not so good \ls\ or a slower but
more powerful one. The answer depends of course on the total amount of
computation time available, and on how the costs improve with
time.

A different type of connection between \ails, \asa\ and \ats\ arises
from certain similarities in the algorithms.  For example, \asa\ can
be seen as an \ails\ without a local search phase (\asa\ samples the
original space $\mathcal{S}$ and not the reduced space
$\mathcal{S}^*$) and where the acceptance criteria is
$\mbox{\ac{LSMC}}(s^*,s^{*}{'},\mbox{\em history\/})$. While \asa\
does not employ memory, the use of memory is the main feature of \ats\
which makes a strong use of historical information at multiple levels.
Given its effectiveness, we expect this kind of approach for
incorporating memory to become widespread in future \ails\
applications.\footnote{In early \ats\ publications, proposals similar
to the use of perturbations were put forward under the name {\em
random shakeup}~\cite{Glo86}. These procedures where characterized as
a ``randomized series of moves that leads the heuristic (away) from
its customary path''~\cite{Glo86}. The relationship to perturbations
in \ails\ is obvious.} Furthermore, \ats, as one prototype of a memory
intensive search procedure, can be a valuable source of inspiration
for deriving \ails\ variants with a more direct usage of memory; this
can lead to a better balance between intensification and
diversification in the search.\footnote{Indeed, in~\cite{Glo89},
Glover uses ``strategic oscillation'' strategies whereby one cycles
over these procedures: the simplest moves are used till there is no
more improvement, and then progressively more advanced moves are
used.}  Similarly, \ats\ strategies may also be improved by features
of \ails\ algorithms and by some insights gained from the research on
\ails.

\subsection{Multi-start based metaheuristics}

Multi-start based metaheuristics can be classified into {\em
constructive\/} metaheuristics and {\em perturbation-based\/}
metaheuristics.

Well-known examples of constructive metaheuristics are \aco\ and
\agrasp\ which both use a probabilistic solution construction
phase. An important difference between \aaco\ and \agrasp\ is that
\aaco\ has an indirect memory of the search process which is used to
bias the construction process, whereas \agrasp\ does not have that
kind of memory. An obvious difference between \ails\ and constructive
metaheuristics is that \ails\ does not construct soutions. However,
both generate a sequence of solutions, and if the constructive
metaheuristic uses an embedded \ls, both go from one local minimum to
another. So it might be said that the perturbation phase of an \ails\
is replaced by a (memory-dependent) construction phase in these
constructive metaheuristics. But another connection can be made:
\ails\ can be used instead of the embedded ``local search'' in an
algorithm like \aco\ or \agrasp. This is one way to generalize \ails,
but it is not specific to these kinds of metaheuristics: whenever one
has an embedded \ls, one can try to replace it by an \ils.

Perturbation-based metaheuristics differ in the techniques they use to
actually perturb solutions. Before going into details, let us
introduce one additional feature for classifying metaheuristics: we
will distinguish between population-based algorithms and those that
use a single current solution (the population is of size 1). For
example, EA, scatter search, and ant colony optimization are
population-based, while \ails\ uses a single solution at each
step. Whether or not a metaheuristics is population-based is important
for the type of perturbation that can be applied. If no population is
used, new solutions are generated by applying perturbations to single
solutions; this is what happens for ILS and VNS.
If a population is
present, one can also use the possibility of recombining several
solutions into a new one. Such combinations of solutions are
implemented by ``crossover'' operators in EAs or in the recombination
of multiple solutions in scatter search.

In general, population-based metaheuristics are more complex to use
than those following a single solution: they require mechanisms to
manage a population of solutions and more importantly it is necessary
to find effective operators for the combination of solutions. Most
often, this last task is a real challenge.  The complexity of these
population-based local search hybrid methods can be justified if they
lead to better performance than non-population based
methods. Therefore, one question of interest is whether using a
population of solutions is really useful.  Unfortunately, there are
very few systematic studies which address this
issue~\cite{CraPoWa98:IJC,GlaPot96,JohMcG97,StuGruLinRut00:ppsn,VaeAarLen96}.
Clearly for some problems such as the TSP with high cost-distance
correlations, the use of a single element in the population leads to
good results, so the advantage of population-based methods is small or
nil.  However, for other problems (with less cost-distance
correlations), it is clear that the use of a population is an
appropriate way to achieve search space diversification. Thus
population based methods are desirable if their complexity is not
overwhelming. Because of this, population-based extensions of \ails\
are promising approaches.

To date, several population-based extensions of \ails\ {\it have} been
proposed~\cite{AppBixChvCoo00,HonKahMoo97,Stu98:phd}. The approaches
proposed in~\cite{HonKahMoo97,Stu98:phd} keep the simplicity of \ails\
algorithms by maintaining unchanged the perturbations: one parent is
perturbed to give one child. But given that there is a population, the
evolution depends on competition among its members and only the
fittests survive.  One can give up this simplicity as was done in the
approach of Applegate et al.~\cite{AppBixChvCoo00}.  Given the
solutions in a population that have been generated by an \ails, they
define a smaller instance by freezing the components that are
in common in all parents. (They do this in the context of the TSP; the
subtours that are in common are then fixed in the sub-problem.)  They
then reoptimize this smaller problem using \ails.  This idea is tested
in~\cite{AppBixChvCoo00}, and they find very high quality solutions,
even for large TSP instances.

Finally, let us discuss variable neighborhood search (VNS) which is
the metaheuristic closest to \ails.  VNS begins by observing that the
concept of local optimality is conditional on the neighborhood
structure used in a local search. Then VNS systemizes the idea of
changing the neighborhood during the search to avoid getting stuck in
poor quality solutions.  Several VNS variants have been proposed. The
most widely used one, {\em basic VNS}, can, in fact, be seen as an
\ails\ algorithm which uses the \ac{Better} acceptance criterion and a
systematic way of varying the perturbation strength. To do so, basic
VNS orders neighborhoods as ${\cal N}_1, \ldots ,{\cal N}_m$ where the
order is chosen according to the neighborhood size. Let $k$ be a
counter variable, $k=1, 2, \ldots ,m$, and initially set $k=1$. If the
perturbation and the subsequent local search lead to a new best
solution, then $k$ is reset to $1$, otherwise $k$ is increased by
one. We refer to ~\cite{HanMla99:mic,MlaHan97:cor} for a description
of other VNS variants.

A major difference between \ails\ and \avns\ is the philosophy
underlying the two metaheuristics: \ails\ explicitly has the goal of
building a walk in the set of locally optimal solutions, while \avns\
algorithms are derived from the idea of systematically changing
neighborhoods during the search.

Clearly, there are major points in common between most of today's high
performance metaheuristics. How can one summarize how \ils\ differs
from the others?  We shall proceed by enumeration as the diversity of
today's metaheuristics seems to forbid any simpler approach. When
comparing to \aaco\ and \agrasp, we see that \ails\ uses perturbations
to create new solutions; this is quite different in principle and in
practice from using construction.  When comparing to EAs and scatter
search, we see that \ails, as we defined it, has a population size of
$1$; therefore no recombination operators need be defined.  We could
continue like this, but we cannot expect the boundaries between all
metaheuristics to be so clear-cut. Not only are hybrid methods very
often the way to go, but most often one can smoothly go from one
metaheuristic to another. In addition, as mentioned at the beginning
of this chapter, the distinction between heuristic and metaheuristic
is rarely unambiguous. So our point of view is not that \ils\ has
essential features that are absent in other metaheuristics; rather,
when considering the basic structure of \ils, some simple yet powerful
ideas transpire, and these can be of use in most metaheuristics, being
close or not in spirit to \ils.

\section{Conclusions}

\ails\ has many of the desirable features of a metaheuristic: it is
simple, easy to implement, robust, and highly effective. The essential
idea of \ails\ lies in focusing the search not on the full space of
solutions but on a smaller subspace defined by the solutions that are
locally optimal for a given optimization engine.  The success of
\ails\ lies in the {\em biased} sampling of this set of local
optima. How effective this approach turns out to be depends mainly on
the choice of the local search, the perturbations, and the acceptance
criterion.  Interestingly, even when using the most na\"\i ve
implementations of these parts, \ails\ can do much better than random
restart. But with further work so that the different modules are well
adapted to the problem at hand, \ails\ can often become a competitive
or even state of the art algorithm. This dichotomy is important
because the optimization of the algorithm can be done progressively,
and so \ails\ can be kept at any desired level of simplicity.  This,
plus the modular nature of \ils, leads to short development times and
gives \ails\ an edge over more complex metaheuristics in the world of
industrial applications. As an example of this, recall that \ails\
essentially treats the embedded heuristic as a black box; then
upgrading an \ails\ to take advantage of a new and better local search
algorithm is nearly immediate.  Because of all these features, we
believe that \ails\ is a promising and powerful algorithm to solve
real complex problems in industry and services, in areas ranging from
finance to production management and logistics.  Finally, let us note
that although all of the present review was given in the context of
tackling combinatorial optimization problems, in reality much of what
we covered can be extended in a straight-forward manner to continuous
optimization problems.

Looking ahead towards future research directions, we expect \ails\ to
be applied to new kinds of problems.  Some challenging examples are:
(i) problems where the constraints are very severe and so most
metaheuristics fail; (ii) multi-objective problems, bringing one
closer to real problems; (iii) dynamic or real-time problems where the
problem data vary during the solution process.

The ideas and results presented in this chapter leave many questions
unanswered. Clearly, more work needs to be done to better understand
the interplay between the \ails\ modules \genis , \genps , \genls ,
and \genps . In particular, we expect significant improvements to
arise through the intelligent use of memory, explicit intensification
and diversification strategies, and greater problem-specific
tuning. The exploration of these issues has barely begun but should
lead to higher performance \ils\ algorithms.

\section*{Acknowledgments}
O.M. acknowledges support from the Institut Universitaire de France.
This work was partially supported by 
the ``Metaheuristics Network'', a
Research Training Network funded by the Improving Human Potential
programme of the CEC, grant HPRN-CT-1999-00106. The information
provided is the sole responsibility of the authors and does not
reflect the Community's opinion.  The Community is not responsible for
any use that might be made of data appearing in this publication.

\newpage

\bibliography{ils} \bibliographystyle{plain}

\end{document}